\documentclass[10pt]{article}  
\usepackage{latexsym,amssymb,amsmath,eepic}

\def\Z{\mathbb{Z}}
\def\R{\mathbb{R}}
\def\C{\mathbb{C}}

\def\proof{\par\medskip\noindent{\em Proof. }}
\def\eproof{\hfill{$\Box$}\bigskip}

\def\de{\delta}

\newtheorem{thm}{Theorem}
\newtheorem{pro}[thm]{Proposition}
\newtheorem{cor}[thm]{Corollary}

\begin{document}
\title{On the theorem of Duminil-Copin and Smirnov about the number of self-avoiding walks in the hexagonal 
lattice} 
\author{Martin Klazar\footnote{Department of Applied Mathematics, Charles University, Faculty of Mathematics and
Physics, Malostransk\'e n\'am\v est\'\i\ 25, 118 00 Praha, Czech Republic, {\tt klazar@kam.mff.cuni.cz}}}
\maketitle

\begin{abstract}
This is an exposition of the theorem from the title, which says that the number of self-avoiding 
walks with $n$ steps in the hexagonal lattice has asymptotics $(2\cos(\pi/8))^{n+o(n)}$. We lift the
key identity to formal level and simplify the part of the proof bounding the growth constant from below.
In our calculation the lower bound comes from an identity asserting that a linear combination of $288$ generating functions 
counting self-avoiding walks in a certain domain by length, final edge direction and winding number modulo $48$
equals the geometric series $\sum_{n\ge1}(2\cos(\pi/8))^nx^n$.
\end{abstract}

\section{Introduction}

Let $s_n$ be the number of $n$-step self-avoiding walks with fixed initial vertex in the hexagonal lattice graph, 
which is a plane graph whose vertices and edges are the vertices and unit 
segments in the tiling of the plane into regular hexagons. In particular, every vertex has degree three. An 
argument using submultiplicativity of $s_n$ proves existence of the limit $\mu=\lim s_n^{1/n}$; clearly, $1\le\mu\le2$.
In 1982, Nienhuis \cite{nien} conjectured that $\mu=\sqrt{2+\sqrt{2}}$. His conjecture was recently proved by 
Duminil-Copin and Smirnov \cite{dumi_smir}. (For more information on counting self-avoiding walks in lattices see 
\cite{gutt_al}.)

\begin{thm}[Duminil-Copin and Smirnov, 2010]\label{dcsthm}
The number $s_n$ of self-avoiding walks in the hexagonal lattice with length $n$ and fixed initial vertex 
has for $n\to\infty$ asymptotics
$$
s_n=\mu^{n+o(n)},\ \mu=2\cos(\pi/8)=\sqrt{2+\sqrt{2}}=1.84775\dots\ .
$$
\end{thm}

\noindent
Here we give an exposition of this remarkable result, the proof of which uses a beautiful combination of geometric 
and generating functions based ideas.

My motivation as an enumerative combinatorialist fond of generating functions 
was to understand the arguments of \cite{dumi_smir} from this angle. Pondering the identities in Lemmas 1 and 2 
of \cite{dumi_smir} one sees that they can be obtained as specializations of a more general formal identity. We 
present this generalization in Proposition~\ref{formiden}. The formal approach treats finite and infinite domains
of vertices on equal footing and skips some limit transitions of \cite{dumi_smir}. The lower 
bound on $\mu$ comes from it in a more direct way compared to \cite{dumi_smir}, by applying the identity 
of Proposition~\ref{formiden} to an appropriate domain. For this domain we get a formal identity stated 
in Corollary~\ref{dusledecek}: If one pigeonholes the length generating function of self-avoiding walks in the domain 
(with fixed starting edge) according to the six possible directions of the final edge and the winding number modulo $48$, then a linear combination of these $288$ generating functions equals the geometric 
series $\sum_{n\ge1}(2\cos(\pi/8))^nx^n$. In deriving the upper bound on $\mu$ in Section 3 we follow the 
arguments of \cite{dumi_smir} more closely. Finally, in Section 4 we prove results on winding number required for the proof, 
such as that the winding number of closed self-avoiding walks attains only two values. 

We adapt the notation of \cite{dumi_smir} to our needs. For example, the hexagonal graph $H$ we work with is not exactly the hexagonal 
lattice graph but its isomorphic embedding in the square lattice graph. We found it more convenient to think in terms 
of rectangular pictures and, for exposing the underlying combinatorial reasoning, 
to divorce the hexagonal lattice graph from its particular embedding in the complex plane. 
This is the reason for using the direction labelling $\ell$
that just produces the complex vectors $v-u$ for oriented edges $uv$ of the hexagonal lattice graph in $\C$.
Also, instead of the self-avoiding walks between mid-edges of \cite{dumi_smir} we work with self-avoiding walks with 
ordinary edges and thus allow coincidence of the final vertex with an earlier vertex of the walk.

\section{A formal identity and the lower bound on $\mu$}

We begin with some definitions.
Let $G=(V,E)$ be a (simple, undirected) graph, typically infinite but with finite degrees. A {\em walk} in $G$ is any sequence of 
vertices and edges $w=(v_1,e_1,v_2,e_2,\dots,e_n,v_{n+1})$ of $G$ such that $n\ge1$, $e_j=v_jv_{j+1}$ and $e_j\ne e_{j+1}$
(we do not allow backtracks); $n=|w|$ is the {\em length of $w$}, $v_1$ is its {\em first (initial)} and 
$v_{n+1}$ {\em last (final) vertex} and similarly for the extremal edges. The {\em reversal} $w^r$ of $w$ is 
$(v_{n+1},e_n,v_n,e_{n-1},\dots,e_1,v_1)$.
$E(w)=\{e_1,\dots,e_n\}$ is the set of its edges and the final edge is denoted $f(w)=e_n$. If $v_1=v_{n+1}$ then $w$ is a {\em closed walk}. If all vertices $v_j$ 
are distinct (thus $w$ is a path) we say that $w$ is a {\em self-avoiding walk} or, for short, a {\em SAW}. A 
{\em closed SAW} is a closed walk in which $v_1=v_{n+1}$ is the only equality of vertices. A {\em weak SAW} is a walk in $G$
such that all vertices $v_j$ are distinct, with the possible exception of the equality $v_{n+1}=v_j$ for some 
$j\le n$. Thus a weak SAW is a SAW or a closed SAW or a concatenation $w_1w_2$ of a SAW $w_1$ and a closed 
SAW $w_2$ such that $w_1\cap w_2$ is exactly the final vertex of $w_1$, which coincides with the initial/final vertex of $w_2$. 

We work with the infinite plane {\em hexagonal graph} 
$$
H=(V,E)\subset\R^2,\ V=\Z\times\Z,
$$ 
arising from the square lattice graph by omitting half of the horizontal edges, those whose left endvertex has
coordinates with different parity. So $H$ is a rectangular drawing of the hexagonal lattice graph $H'$; it
is a cubic graph whose all faces are rectangular hexagons. We obtain $H$ from $H'$ by keeping the  
unit length of edges and the horizontal edges but compressing $H'$ horizontally so that all non-horizontal 
edges become vertical. For any vertex $v$ we denote by $N(v)=\{e_1,e_2,e_3\}$ the triple of edges incident with $v$ 
with some fixed counter-clockwise indexing (each $N(v)$ has three of them). By $s_n$ we denote the number of SAWs in $H$ with length $n$ and first vertex equal to $(0,0)$. As we mentioned, $\mu=\lim s_n^{1/n}$ exists. 
We present the proof of \cite{dumi_smir} that $\mu=2\cos(\pi/8)$.

Let $X\subset V$ be a finite or infinite connected set of vertices, we call it a {\em domain}.
$H\backslash X$ denotes the graph induced by $H$ on the 
vertex set $V\backslash X$. We call $X$ a {\em simple domain} if every connected component of $H\backslash X$ is infinite. 
$E(X)\subset E$ is the set of
edges with an endvertex in $X$ and $\partial X\subset E(X)$ are the {\em border edges}, the edges with exactly 
one endvertex in $X$. A {\em cycle $C$ in $X$} is any $2$-regular 
connected subgraph of $H$ with all vertices in $X$; $|C|$ is the number of its edges (or vertices).
For a domain $X$ and two edges $a\in\partial X$ and 
$e\in E(X)$ we define $X(a,e)$ to be the set of all weak SAWs $w=(v_1,e_1,\dots,e_n,v_{n+1})$ such that 
$e_1=a$, $v_1$ is the vertex of $a$ not in $X$, $e_n=e$ and if $v_j\not\in X$ then $j=1$ or $j=n+1$. 
Thus $X(a,e)$ is the set of weak SAWs going in $X$ from $a$ to $e$. If $e\in\partial X$ and the vertices of $a$ and $e$ 
in $V\backslash X$ are distinct then $X(a,e)$ consists only of SAWs; if they coincide then $X(a,e)$ consists only 
of closed SAWs. We are interested in the generating function (i.e., formal power series)
$$
F_{a,e}(x)=F^X_{a,e}(x)=\sum_{w\in X(a,e)} x^{|w|}\;.
$$       
Note that if $a=e$ then $F_{a,a}(x)=x$. For any set of edges $A\subset E(X)$ we define the more general set of 
weak SAWs $X(a,A)=\bigcup_{e\in A}X(a,e)$ and the corresponding generating function 
$$
F_{a,A}(x)=\sum_{w\in X(a,A)} x^{|w|}=\sum_{e\in A}F_{a,e}(x)\;.
$$ 

Consider the primitive sixth root of unity
$$
\zeta=\exp(2\pi i/6)
$$
and two mappings
$$
\ell:\;\overline{E}\to\{\zeta^j\;|\;j=0,1,2,3,4,5\}\ \mbox{ and }\ r:\;\mbox{walks in $H$}\to\Z
$$
which we define in detail after stating the identity. $\overline{E}$ means oriented edges (marked by bars), every $e\in E$ appears in
$\overline{E}$ in its two orientations, and we extend $\ell$ to walks by setting 
$\ell(w)=\ell(\overline{e})$ where $e=f(w)$ and has certain orientation. The mapping $r$ 
gives winding numbers of walks. Now we state the main identity and afterwards give further explanation.
\begin{pro}\label{formiden}
Let $H=(V,E)$ be the hexagonal graph, $X\subset V$ a simple domain and $a\in\partial X$ its border edge.
Then we have in $\C[y,y^{-1}][[x]]$ the formal identity
\begin{eqnarray*}
&&\sum_{e\in\partial X}\ell(\overline{e})y^{r(e)}F_{a,e}(x)=
(1+x\zeta^{-2}y+x\zeta^2y^{-1})\sum_{w\in S}\ell(w)y^{r(w)}x^{|w|}+\\
&&+\;(\zeta^2 y^{-4}+\zeta^{-2}y^4)\sum_{w\in S,\;C\in C(w)}\ell(w)y^{r(w)}x^{|w|+|C|}
\end{eqnarray*}
where $S$ is the set of all SAWs that start in $a$ and have all vertices but the first in $X$ and 
$C(w)$ is the set of all cycles in $X$ that intersect $w$ only in its final vertex.
\end{pro}
Here $x$ and $y$ are formal variables and only later we plug in for them some complex numbers. 

We define orientation of edges and their direction $\ell$. 
For $e\in E$, an {\em oriented edge} is a pair $\overline{e}=(e,v)$ where $v$ is an endvertex
of $e$. We think of $\overline{e}$ as oriented away from $v$. The {\em basic orientation} is the left-right and 
bottom-up orientation of edges in $E$ and the {\em basic direction} $\ell_b$ is $\ell_b(e)=\zeta^0=1$ for every horizontal 
edge and $\ell_b(e)=\zeta^1$ (resp. $\ell_b(e)=\zeta^2$) if $e$ is vertical and its top vertex is incident with 
a horizontal edge going to the right (resp. to the left). For general oriented edges we define $\ell$ by 
$\ell(\overline{e})=\ell(e,v)=\ell_b(e)$ if $\overline{e}$ is oriented as in the basic orientation and 
$\ell(\overline{e})=\zeta^3\ell_b(e)=-\ell_b(e)$ if it is oriented in the opposite way. The key features of $\ell$ are that
$$
\ell(e,u)+\ell(e,v)=0
$$
for any edge $e=uv$ and the multiplicative property 
$$
\ell(\overline{e_{i+1}})=\zeta^2\ell(\overline{e_i}),\ i=1,2,3\ (\mbox{mod $3$}),
$$
holding for any triple $N(v)=\{e_1,e_2,e_3\}$ if the $e_i$ are oriented in the same way, towards $v$ or away from it. 
For $e\in\partial X$ we define $\overline{e}=(e,e\cap X)$, so $e$ is oriented away from its vertex in $X$ 
(this concerns the left 
side of the identity). Any walk $w=(v_1,e_1,\dots,e_n,v_{n+1})$ has a natural orientation from $v_1$ to $v_{n+1}$
and we set $\ell(w)=\ell(\overline{e_n})$ where the final edge is oriented in the opposite way to the natural 
orientation of $w$ (this concerns summands on the right side of the identity). 

We define the winding number $r$. A {\em hook} $h=(v_1,e_1,v_2,e_2,v_3)=e_1e_2$ 
is any SAW in $H$ with length $2$; $v_2$ is its {\em central vertex}. We call $h$ a {\em left hook} if it `turns left', $e_2$ 
comes clockwise after $e_1$ with respect to $v_2$; else, if $e_2$ comes counter-clockwise after $e_1$, we call $h$ 
a {\em right hook}. To the left hooks we give weight $r(h)=+1$ and to the right ones $r(h)=-1$.
For any walk $w=(v_1,e_1,\dots,e_n,v_{n+1})$ we define its {\em winding number} by
$$
r(w)=\sum_{i=1}^{n-1}r((v_i,e_i,v_{i+1},e_{i+1},v_{i+2}))\;.
$$
Thus $r(w)$ is the difference between the number of left and right hooks contained in $w$. When $|w|=1$ we set $r(w)=0$.
It follows that $\zeta^{r(w)}=\ell(\overline{e_n})/\ell(\overline{e_1})$, with the natural orientation of the walk. 
We will show that if $X$ is a simple domain and the edges $a,e\in\partial X$ are fixed then all SAWs $w\in X(a,e)$ have the same winding number, which explains the notation $r(e)=r(w)$ on the left side of the identity.

Our explanation of the identity is now complete. Before proving it we review a few properties of
winding number used in the proof. Clearly, $r(w^r)=-r(w)$ and if $w$ is split by a vertex $v$ in two walks as $w=w_1w_2$ 
then
$$
r(w)=r(w_1)+r(w_2)+r(h)
$$
where $h$ is the hook in $w$ with central vertex $v$. If $w=(v_1,e_1,\dots,e_n,v_{n+1}=v_1)$ is a closed walk, it 
is convenient to define its {\em complete winding number $r^*$} as 
$$
r^*(w)=r(w')=r(w)+r(h)
$$
where $w'$ is the extension of $w$ by the edge $e_{n+1}=e_1$ and $h$ is the hook $(v_n,e_n,v_1,e_1,v_2)$. Thus 
$r^*(w)$ is the sum of weights of the hooks in $w$ when their central vertex runs through {\em all} vertices 
of $w$ including the initial/final vertex (the hook of which, strictly speaking, is not contained in $w$). 
It is easy to see that for any closed walk $r^*(w)=6k$ for 
some $k\in\Z$. An important case is when $w$ is a closed SAW. We recall a particular case of 
Jordan's curve theorem.

\begin{pro}\label{jordan}
Every closed SAW $w$ in $H$ separates the plane in two connected components, 
the bounded {\em interior of $w$} and the unbounded {\em exterior of $w$}. For any straight segment $uv$ that intersects 
$w$ in a point different from $u$, $v$ and the vertices of $w$, the endpoints $u$ and $v$ lie in distinct components of 
$\R^2\backslash w$. 
\end{pro}
We say that a closed SAW $w$ is {\em clockwise (resp. counter-clockwise) oriented} if its interior lies to the right (resp. left) of its 
edges when they are naturally oriented.

\begin{pro}\label{konstrot} The following holds.
\begin{enumerate}
\item If $w$ is a closed SAW in $H$ then, depending on the counter-clockwise or clockwise orientation,
$$
r^*(w)=+6\ \mbox{ or }\ r^*(w)=-6\;.
$$
\item If $X\subset V$ is a simple domain and $a,e\in\partial X$ are two boundary edges then all SAWs 
of $X(a,e)$ have the same winding number.
\end{enumerate}
\end{pro}

\noindent
We prove Proposition~\ref{konstrot} in Section 4. For brevity we do not prove Proposition~\ref{jordan} and simply 
assume it, although we could prove it by the same inductive argument of Section 4. 
For more information on Jordan's theorem and its proofs see \cite{hale} and \cite{hale1}. Now we prove the identity.

\bigskip
{\bf Proof of Proposition~\ref{formiden}.} Let $H$, $X$ and $a$ be as given. Let 
$X(a)=X(a,E(X))$, the set of all weak SAWs whose first edge is $a$, first vertex is the endvertex of 
$a$ not in $X$ and all vertices but the initial one and possibly the final one lie in $X$. Consider the formal sum 
$$
\sum_{(w,v)\in X(a)\times X,\;v\in f(w)}\ell(f(w),v)y^{r(w)}x^{|w|}\;.
$$
For each $n$ there are only finitely many pairs $(w,v)$ in the sum with 
$|w|=n$. Hence the sum is a well defined element of $\C[y,y^{-1}][[x]]$ since for any ordering of the set 
of pairs into a sequence $p_1,p_2,\dots$, the corresponding summands $s_1,s_2,\dots$ yield partial sums 
$s_1,s_1+s_2,\dots$ formally converging to the same power series, namely the one with coefficient of $x^n$ equal to
$\sum\ell(f(w),v)y^{r(w)}$, summed over the mentioned finitely many pairs. Moreover, it follows that 
for any partition of $p_1,p_2,\dots$ into subsequences $p_{1,j},p_{2,j},\dots$, $j=1,2,\dots$, we have 
the double counting identity $\sum_i s_i=\sum_i\sum_j s_{i,j}=\sum_j\sum_i s_{i,j}$. The identity in question is such 
double counting identity. On the one hand, grouping together the pairs $(w,v)$ sharing the same weak SAW $w$ we get
the left hand side because for each $w$ with $f(w)=uu'\in E(X)\backslash\partial X$ we have two contributions that cancel 
out ($\ell(uu',u)+\ell(uu',u')=0$) and for each $w$ with $f(w)=e\in\partial X$ we have a single contribution in which 
$r(w)$ depends only on $e$ (by 2 of Proposition~\ref{konstrot}). 

On the other hand, grouping together the pairs $(w,v)$ sharing the same vertex $v$ gives the right hand side of the 
identity. We fix $v\in X$ and split the group of the corresponding weak SAWs as 
$$
X(a,N(v))=T_1\cup T_2
$$  
where $T_1$ consists of the $w$ containing one or two edges of $N(v)$ and $T_2$ of the $w$ containing all three of them. 
Let $N(v)=\{e_1,e_2,e_3\}$ with a fixed counter-clockwise indexing.

Suppose that $w\in T_1$. If $E(w)\cap N(v)$ is $f(w)$ and the preceding edge, by deleting $f(w)$ we get another element 
of $T_1$. If $E(w)\cap N(v)$ is just $f(w)=e_i$, we can always extend $w$ by $e_{i+1}$ and $e_{i+2}$ of $N(v)$ and 
get an element of $T_1$ (now it is important to have weak SAWs and not just SAWs). Thus $T_1$ splits into
triples $w,w',w''$ such that $E(w)\cap N(v)=\{e_i\}=\{f(w)\}$, so $w$ is a SAW and $v$ is its final vertex, 
and $w'$ (resp. $w''$) is an extension 
of $w$ by $e_{i+1}$ (resp. by $e_{i+2}$). If $s=\ell(f(w),v)y^{r(w)}x^{|w|}$ is the summand corresponding to $w$, 
then due to additivity of winding and length and the multiplicative property of $\ell$ the summand of $w'$ 
(resp. $w''$) is $s'=\zeta^2y^{-1}xs$ (resp. $s''=\zeta^4yxs$). Taking $s$ out and summing over all $w$ and $v\in X$, we obtain the first term on the right side. 

Suppose that $w\in T_2$. Let $e_i$ be the edge of $N(v)$ appearing in $w$ first. 
Then $w$ splits in $w=w_1w_2$ where $w_1$ is the initial part with final edge $e_i$ and $w_2$ is a closed SAW with 
initial edge $e_{i+1}$ and final edge $e_{i+2}$ or vice versa. Reverting the direction of $w_2$ we get another element
of $T_2$. This operation is an involution that splits $T_2$ into pairs $w,w'$ where $w$ (resp. $w'$) is obtained from 
a SAW $w_1$ with final vertex $v$ and a cycle $C\in C(w_1)$ by setting $v$ to be the initial/final vertex of $C$ and choosing
one of the two directions, which turns $C$ into a closed SAW $w_2$ (resp. $w_2'$). We may assume that the hook of $w$ 
(resp. $w'$) with the central vertex $v$ is a left one (resp. a right one). By 1 of Proposition~\ref{konstrot}, 
$r(w)=r(w_1)+r(w_2)+1=r(w_1)+(r^*(w_2)-(-1))+1=r(w_1)-4$ and, similarly, $r(w')=r(w_1)+4$ because $w_2$ is 
clockwise oriented (note that $a$ and hence $e_i$ lie in the exterior of $w_2$). We see that if $t$ is the summand 
of $w_1$ then the summand of $w$ (resp. $w'$) is $s=\zeta^2y^{-4}x^{|C|}t$ (resp. $s'=\zeta^4y^4x^{|C|}t$). 
Taking $x^{|C|}t$ out and summing over all $w_1,C$ and $v\in X$, we obtain the second term on the right side.
\eproof

Before we apply the identity to concrete domains, we point out two simplifications of its left side and discuss them
in some generality. The first uses that $r(e)$ often depends only on the direction $\ell(\overline{e})$. 
We introduce some terminology. A biinfinite path $p$ in $H$ is a {\em quasiline} if left and right hooks alternate on it; 
such $p$ separates $H$ in two components. A SAW $w$ is a {\em quasisegment} if $|w|\ge 2$ and is contained in a quasiline. 
A {\em quasihalfline} is defined in an obvious way. A set of vertices is {\em straight} if it forms a quasisegment, a 
quasihalfline or a quasiline. Any straight set of vertices has a unique extension to a quasiline. Suppose that 
$X\subset V$ is a domain and $a\in\partial X$. We say that a partition of border edges
$$
\partial X\backslash\{a\}=A_1\cup A_2\cup\dots\cup A_k
$$
is a {\em straight partition} if there exist vertex sets $U_i\subset X$, $i=1,2,\dots,k$, such that, for each $i$, 
$U_i$ is straight and the edges in $A_i$ are incident with $U_i$, lie all on one side of the quasiline extending $U_i$
(and not on it) and form an interval (if $e,f\in A_i$ with endvertices $u,v\in U_i$ and $g\in E(X)$ is incident to $U_i$,
lies on the same side of $U_i$ as $e$ and $f$ and its vertex in $U_i$ is between $u$ and $v$, then $g\in A_i$).

\begin{pro}\label{zjednodus}
Let $X\subset V$ be a simple domain, $a\in\partial X$ and $\partial X\backslash\{a\}=A_1\cup A_2\cup\dots\cup A_k$ be 
a straight partition. Then $e,f\in A_i$ implies that $\ell(\overline{e})=\ell(\overline{f})$ and $r(e)=r(f)$, that is, 
$r(w)=r(w')$ whenever $w\in X(a,e)$ and $w'\in X(a,f)$. Thus the left side of the identity can be written as 
$$
\ell(\overline{a})x+\sum_{i=1}^k\ell(A_i)y^{r(A_i)}F_{a,A_i}(x)
$$
where $\ell(A_i)=\ell(\overline{e})$ and $r(A_i)=r(e)$ for any representing edge $e\in A_i$.
\end{pro}

\noindent
The proof is in Section 4. Further simplification uses the mirror symmetry along the $x$-axis. It is an automorphism of 
$H$ and we denote it by $M$. Suppose that $X$ is a simple domain such that $M(X)=X$, $a=(0,0)(1,0)\in\partial X$ and 
$\partial X\backslash\{a\}=A_1\cup A_2\cup\dots\cup A_k$ is a straight partition. Such partition can be usually defined 
so that all or almost all parts come in symmetric pairs. In any case, 
if $A$ and $B=M(A)$, $A\ne B$, are two parts of such symmetric pair, we can pair their summands on the left side 
of the identity and get
$$
\ell(A)y^{r(A)}F_{a,A}(x)+\ell(B)y^{r(B)}F_{a,B}(x)=
\frac{\ell(A)y^{r(A)}+(\ell(A)y^{r(A)})^{-1}}{2}F_{a,A\cup B}(x)\;.
$$
This is because $M$ changes signs of winding numbers of SAWs and conjugates directions of border edges 
(in orientation away from $X$) and gives a length-preserving bijection between $X(a,e)$ and $M(X)(M(a),M(e))=X(a,M(e))$ 
for any $e\in\partial X$. 

Let us proceed with the lower bound on $\mu$. If we specialize $y$ to a complex number satisfying 
$\zeta^2y^{-4}+\zeta^{-2}y^4=0$, the second term on the right side of the identity in Proposition~\ref{formiden} 
disappears. This equation has eight solutions
$$
y_j=\zeta^{1/8}=\exp(2\pi i(6j+1)/48),\ 0\le j\le 7\;.
$$
For exactly the same reason we can make almost vanish the left side of the identity by selecting appropriate $X$. We set 
$$
X=\Z^2\backslash\{(p,q)\;|\;p\le0,\; p-1\le q\le1-p\}\ \mbox{ and }\ a=(0,0)(1,0)\in\partial X\;.
$$
We call $X$ the {\em slit plane domain}. It is a simple domain, symmetric along the $x$-axis and its border edges, which 
are all except $a$ vertical, have
straight partition $\partial X\backslash\{a\}=A_1\cup A_2$ where $A_1$ are the ones above the $x$-axis and $A_2=M(A_1)$ are the symmetric ones. The corresponding set $U_1\subset X$ is the quasihalfline with vertices 
$(0,2),(0,3),(-1,3),(-1,4),(-2,4),(-2,5),\dots,$ and $U_2=M(U_1)$ is the symmetric quasihalfline. It is easy to see that 
$\ell(A_1)=\zeta^4$ and $r(A_1)=4$. Also, $\ell(\overline{a})=-1$. 
By Proposition~\ref{zjednodus}, the symmetry principle and the equation defining $y_j$, for each $y=y_j$ the left side 
of the identity for the slit plane domain equals
$$
\sum_{e\in\partial X}\ell(\overline{e})y_j^{r(e)}F_{a,e}(x)=\ell(\overline{a})x+\frac{(\zeta y_j)^4+
(\zeta y_j)^{-4}}{2}F_{a,A_1\cup A_2}(x)=-x\;.
$$
\begin{cor}\label{dusledecek}
Let $X\subset V$ be the slit plane domain and $a=(0,0)(1,0)\in\partial X$. For each $l\in\{0,1,\dots,5\}$ and 
$k\in\{0,1,\dots,47\}$ we consider the generating function
$$
G_{l,k}(x)=\sum_{\ell(w)=\zeta^l,\;r(w)\equiv k\;({\rm mod}\;48)}x^{|w|}
$$
where we sum over the SAWs $w$ starting in $a$ and with all vertices but the first in $X$. Then for each
$j\in\{0,1,\dots,7\}$ we have in $\C[[x]]$ the formal identity
$$
\frac{x}{1+2\mathrm{Re}(\zeta^4y_j)x}=-\sum_{l=0}^5\sum_{k=0}^{47}\zeta^ly_j^kG_{l,k}(x)\;.
$$
For $j=7$ it follows that $G_{l,k}(x_c)=+\infty$ for some $l,k$ and $x_c=\frac{1}{2\cos(\pi/8)}$. 
\end{cor}
\proof
We SAW that for $y=y_j$ and the slit plane domain $X$ the identity of Proposition~\ref{formiden} turns in 
$$
-x=(1+\zeta^{-4}xy_j^{-1}+\zeta^4xy_j)\sum_{w\in S}\ell(w)y_j^{r(w)}x^{|w|}\;.
$$
Rearranging the equation and using that each $y_j$ is a (primitive) $48$-th root of unity we get the stated formal 
identity. For $j=7$ the real part of $\zeta^4y_j=\exp(2\pi i(2/3+43/48))=\exp(9\pi i/8)$ equals $-\cos(\pi/8)$.
For $j=7$ and $x=x_c=(2\cos(\pi/8))^{-1}$ the geometric series on the left side  diverges, 
which implies that some of the $288$ power series $G_{l,k}(x)$ on the right side diverges as well. 
\eproof

\noindent
Thus $G_{l,k}(x)=\sum_{n\ge1}a_nx^n$ with this special $k$ and $l$ satisfies $\limsup a_n^{1/n}\ge2\cos(\pi/8)$ and so, 
since $s_n\ge a_n$, $\mu\ge2\cos(\pi/8)$.

\section{The upper bound on $\mu$}

We show that $s_n$, the number of $n$-step SAWs in $H$ starting at the origin, has for $n=1,2,\dots$ and every 
$\de>0$ bound 
$$
s_n\ll_{\de}(x_c^{-1}+\de)^n=(2\cos(\pi/8)+\de)^n\;.
$$
This implies that $\mu\le2\cos(\pi/8)$.

We apply Proposition~\ref{formiden} to finite trapezoidal domains that we now define. For integers $r,s\ge2$ we set 
$$
X_{r,s}=\{(p,q)\in\Z^2\;|\;1\le p\le r,\;|q|\le 2s+p-2\}\ \mbox{ and }\ a=(0,0)(1,0)\in\partial X_{r,s}\;. 
$$
$X=X_{r,s}$ is a simple domain, symmetric along the $x$-axis and its border edges have straight partition
$$
\partial X\backslash\{a\}=A^+\cup B^+\cup C\cup B^-\cup A^-
$$
where $A^+$ are the horizontal edges above $a$, $B^+$ are the vertical edges above the $x$-axis, $C$ are the 
horizontal edges to the right of $X$ and $B^-=M(B^+)$ and $A^-=M(A^+)$ are the symmetric parts. The corresponding
five quasisegments in $X$ are, in the respective order,
$U_1=\{(1,q)\;|\;1\le q\le 2s-1\}$, $U_2=\{(p,2s+p-2),(p,2s+p-3)\;|\;1\le p\le r\}$, $U_3=\{(r,q)\;|\;|q|\le 2s+r-2\}$ 
and the symmetric ones $U_4=M(U_2)$ and $U_5=M(U_1)$. Clearly, $\ell(\overline{a})=-1$, $\ell(A^+)y^{r(A^+)}=-y^3$, $\ell(B^+)y^{r(B^+)}=\zeta^2y^2$ and $\ell(C)y^{r(C)}=1$. Denoting $A=A^+\cup A^-$ and $B=B^+\cup B^-$, we get the 
next result, Lemma 2 of \cite{dumi_smir}. 

\begin{cor}\label{finidsq}
Let $x_c=\frac{1}{2\cos(\pi/8)}$. For all integers $r,s\ge2$ we have for the domain $X=X_{r,s}$ with border edges
partitioned into $\{a\}\cup A\cup B\cup C$ the identity
$$
\cos(3\pi/8)F_{a,A}(x_c)+\cos(\pi/4)F_{a,B}(x_c)+F_{a,C}(x_c)=x_c\;.
$$
\end{cor}
\proof
In Proposition~\ref{formiden} we set $X=X_{r,s}$, $y=y_7$ and $x=x_c$, which 
is admissible since all sums are finite. Then both factors on the right side of the identity are zero and the right side vanishes completely. By Proposition~\ref{zjednodus} and by symmetry we obtain equation
$$
-x_c-\frac{y_7^3+y_7^{-3}}{2}F_{a,A}(x_c)+\frac{\zeta^2y_7^2+\zeta^{-2}y_7^{-2}}{2}F_{a,B}(x_c)+F_{a,C}(x_c)=0\;.
$$
The coefficient of $F_{a,A}(x_c)$ equals $\cos(3\cdot2\pi(43/48))=-\cos(3\pi/8)$ and that of $F_{a,B}(x_c)$ equals 
$\cos(2\cdot2\pi(1/6+43/48))=\cos(\pi/4)$. Rearrangement gives the stated formula.
\eproof

\noindent
It only is important that $\cos(3\pi/8),\cos(\pi/4)$ and $x_c$ are positive and the involved power 
series have nonnegative coefficients. This identity will provide the crucial estimate.

We define the set of border edges $D\subset\partial X_{r,s}$ as $D=B\cup C\cup\{e,f\}$ where $e$ is the top 
edge of $A^+$ and $f$ the bottom edge of $A^-$. Note that $D$ are exactly the border edges incident with the set 
$Y\subset X_{r,s}$, defined as $Y=U_2\cup U_3\cup U_4$. We call $U_1\cup\{(1,0)\}\cup U_5$ the {\em left side}, 
$U_2$ the {\em top side}, $U_3$ the {\em right side} 
and $U_4$ the {\em bottom side of $X_{r,s}$}. Note that each SAW starting at $a$ needs at least $s$ steps to reach the 
top or bottom side of $X_{r,s}$ and at least $r$ steps to reach the right side. We have this stratification property of trapezoids: If $X=X_{r,s}$, $X'=X_{r-1,s}$ and $X''=X_{r,s-1}$, then $X',X''\subset X$, $X\backslash X'$ 
is the right side of $X$ and $X\backslash X''$ is the union of the top and bottom side of $X$.

Let $a=(0,0)(1,0)$ be the initial edge as before. We bound $s_n$ from above by considering four classes of SAWs in 
$H$ with decreasing degree of generality and by counting their SAWs by length $n$. 
All SAWs starting at $(0,0)$ are counted by $s_n$.
The {\em SAWs starting at $a$} start at $(0,0)$, have initial edge $a$ and are counted by $b_n$.
{\em Halfplane SAWs}, counted by $c_n$, are those of them that after the initial edge $a$ stay in the right 
halfplane $x\ge1$. {\em Trapezoidal SAWs}, counted by $d_n$, is the further subset
$$
\bigcup_{r,s=2}^{\infty}X_{r,s}(a,D)
$$ 
where the trapezoids $X_{r,s}$ and sets of border edges $D=D_{r,s}$ were defined above. For example, $d_1=d_2=d_3=0$ and 
$d_4=2$. We estimate the numbers $s_n,b_n,c_n$ and $d_n$.

For this we employ and first recall the symmetry properties of $H$. For every two edges 
$e,f\in E$ and their endvertices $u\in e,v\in f$ there exist exactly two automorphisms of $H$ sending $u$ to $v$ and 
$e$ to $f$, one preserving the counter-clockwise order of edges in $N(t)$ for any vertex $t$ and the other reversing 
it. The former automorphisms are of the {\em first kind} and the latter of the {\em second kind}. The inverse to 
an automorphism is clearly of the same kind. For example, the mirror symmetry along the $x$-axis $M$ is of the 
second kind and is an involution. In the original hexagonal 
lattice graph $H'$ automorphisms of the first kind are induced by shifting the plane to identify 
$u$ with $v$ and then rotating it around $v$ to identify the shifted $e$ with $f$; automorphisms 
of the second kind need additional mirror symmetry along 
the line going through $f$. If $p$ is a quasiline in $H$ and $e$ is an edge incident with $p$ but not on $p$, we 
denote by $H(p,e)$ the set of SAWs $w$ in $H$ starting at $e$ and such that $p$ separates $e$ and the remaining 
edges of $w$. Clearly, both automorphisms sending $v=e\cap p$ to the origin and $e$ 
to $a$ provide a length-preserving bijection between $H(p,e)$ and the set of halfplane SAWs defined above.
As we know, the border edges in $A\cup\{a\}$, $B^+$, $C$ and $B^-$ are separated from $X_{r,s}$ by the quasiline 
extending, respectively, the left side, the top side, the right side and the bottom side of $X_{r,s}$. Thus $H(p,e)$, 
for $p$ extending a side of $X_{r,s}$ and $e$ in the corresponding set of border edges, is in length-preserving 
bijection with halfplane SAWs. Also, every SAW starting in such $e$ and with all vertices but the first in $X_{r,s}$
lies in the set $H(p,e$). 

Let us relate the numbers $s_n,b_n,c_n$ and $d_n$. We have the simple bounds
$$
s_n=3b_n\ \mbox{ and }\ b_n\le c_{n+1}+\sum_{i=1}^{n+3} c_ic_{n+4-i},\ n=1,2,\dots\;.
$$
The equality follows from the fact that every edge incident with the origin is sent by an automorphism 
fixing the origin to $a$. As for the inequality, for a SAW $w$, $|w|=n\ge1$, starting at $a$ consider the rightmost 
vertical quasiline $p$ containing its vertex and take the topmost of such vertices $v$. We distinguish two cases 
depending on whether $v$ is incident to an edge $e$ to the right of $p$ or not. If $v$ is incident to such edge $e$
then $v$ is the last vertex of $w$ for else the hook of $w$ with the central vertex $v$ would be in $p$ and $v$ would
not be the topmost vertex of $w$ on $p$. We extend $w$ by the one-step SAW $w_0=e$. In the other case we step 
from $v$ first up and then to the right of $p$ (this does not revisit $w$) and let $w_0$ be the corresponding two-step SAW. 
In the first case we set $w_1=F((ww_0)^r)$ and in the second $w_2=F((w'w_0)^r)$ and $w_3=F(w_0^rw'')$ where $w'$ is the initial part of $w$ ending in $v$, $w''$ 
(which may be empty) is the final part of $w$ starting in $v$ and $F$ is the unique automorphism of $H$, say of the 
first kind, that sends the final vertex of $w_0$ to the origin and $f(w_0)$ to $a$. All three SAWs $w_1,w_2$ and 
$w_3$ are halfplane SAWs because $(ww_0)^r,(w'w_0)^r$ and $w_0^rw''$ are in $H(p,f(w_0))$. Thus $w\mapsto w_1$ or 
$w\mapsto(w_2,w_3)$ is an injective mapping sending $w$ to a halfplane SAW or to a pair of them and $|w_1|=n+1$ or  
$|w_2|+|w_3|=n+2|w_0|=n+4$. It is injective because if $w_1$ or $(w_2,w_3)$ is in the image of the mapping, 
$w$ is reconstructed by applying the unique automorphism $F'$ of the first kind sending the origin to the last vertex of 
$w_1$ or $w_2$ and $a$ to $f(w_1)$ or $f(w_2)$ and by omitting the last edge of $F'(w_1)$ or the two edges of 
$F'(w_2)\cap F'(w_3)$. This proves the inequality.  

The breakthrough bound
$$
d_n<4nx_c^{-n},\ n=1,2,\dots,
$$
follows from Corollary~\ref{finidsq}. Consider a trapezoidal SAW $w\in X_{r,s}(a,D)$ with length $n$. If 
$f(w)\in A\cup B$ then, by definition of $X_{r,s}$, $w\in X_{n,s}(a,D)$ and if $f(w)\in C$ then 
$w\in X_{r,n}(a,D)$. Thus all trapezoidal SAWs with length $n$ are contained in the union of $2n-2$ 
sets
$$
\bigcup_{r=2}^n X_{r,n}(a,D)\cup\bigcup_{s=2}^n X_{n,s}(a,D)\;.
$$
For fixed $r$ and $s$, the number of SAWs $w\in X_{r,s}(a,D)$ with length $n$ is less than $x_c^{-n}$ because the identity 
in Corollary~\ref{finidsq} implies that $F_{a,D}(x_c)\le F_{a,A}(x_c)+F_{a,B}(x_c)+F_{a,C}(x_c)\le x_c/\cos(3\pi/8)<2$ 
and the coefficient of $x^n$ in $F_{a,D}(x)$ cannot be $\ge 2x_c^{-n}$. Summing over the $2n-2$ sets of SAWs 
gives the bound.

To close the chain we bound $c_n$ in terms of $d_n$, halfplane SAWs in terms of trapezoidal SAWs.
Arguing as before we derive the recurrent inequality
$$
c_n\le \sum_{i=\sqrt{n/6}}^{n+1} d_ic_{n+2-i}+\sum_{i=\sqrt{n/6}}^{n+3} d_ic_{n+4-i},\ n=3,4,\dots\;.
$$
Let $w$ be a halfplane SAW with $|w|=n\ge3$. Consider the trapezoid $X_{r,s}$ containing $w$  
(except the initial vertex) that has the minimum size $|X_{r,s}|$; we have $|X_{r,s}|\le r(4s+2r-3)<6\max(r,s)^2$. 
Due to the minimality of $X_{r,s}$, $|w|\ge3$ and the stratification property of trapezoids, $w$ visits the right 
side of $X_{r,s}$ or the union of its top and bottom side. If only the former visit occurs then $s=2$ and 
if only the latter visit occurs then $r=2$. In any case, if $v_j$ is the last vertex of $w$ in 
$Y=Y_{r,s}$ then $j\ge\max(r,s)$ because $w$ needs so many steps to reach the right side 
and the top or bottom side. Also, $\max(r,s)\ge\sqrt{n/6}$ for else $|X_{r,s}|<n$ and $w$ would not fit in 
$X_{r,s}$. Thus $j\ge\sqrt{n/6}$. As above we distinguish two cases depending on whether 
$v_j$ is incident with a border edge in $D$ or not. If $v_j$ is incident with such edge $e\in D$, we set 
$w_0=e$ (oriented away from $X_{r,s}$), $w_1=w'w_0$ and $w_2=F(w_0^rw'')$ where $w'$ is the initial part of $w$ ending 
in $v_j$, $w''$ (which may be empty) is the final part of $w$ starting in $v_j$ and $F$ is the unique automorphism of 
$H$, say of the first kind, that sends the final vertex of $w_0$ to the origin and $f(w_0)=e$ to $a$. (Unlike the bound for $b_n$ and $c_n$, 
here if $e\in A$ then $v_j$ need not be the last vertex of $w$.) If $v_j$ is not incident with a border edge in $D$, it has (at 
least) two neighbors in $Y_{r,s}$ and each of them is incident with an edge in $D$. It follows that one of these 
neighbors is not 
in $w$, for else both would be earlier vertices of $w$ and both hooks in $w$ with these central vertices would use edges
incident with $v_j$, both of which would be in the natural orientation oriented towards $v_j$, which is impossible. Hence there exists a two-step SAW $w_0$ starting in $v_j$, not revisiting $w$ and with final edge in $D$. 
We define $w_3=w'w_0$ and $w_4=F(w_0^rw'')$ where $w',w''$ and $F$ are as in the first case. Clearly, $|w_1|+|w_2|=n+2$, 
$|w_3|+|w_4|=n+4$, $w_1$ and $w_3$ are trapezoidal SAWs, $|w_1|,|w_3|\ge j\ge\sqrt{n/6}$, and $w_2$ and $w_4$ are halfplane SAWs because $w_0^rw''$ is 
in $H(p,f(w_0))$ where $p$ is the quasiline containing a side of $X_{r,s}$ and separating $X_{r,s}$ from the last vertex of 
$w_0$. The mapping $w\mapsto (w_1,w_2)$ or $w\mapsto (w_3,w_4)$ is injective because one inverts it by applying
the unique automorphism of the first kind $F'$ sending the origin to the last vertex of $w_2$ or $w_4$ and $a$ to 
$f(w_2)$ or $f(w_4)$ and by deleting the one or two edges of $F'(w_2)\cap w_1$ or $F'(w_4)\cap w_3$. This proves the inequality.

Combining the bounds we prove the upper bound on $s_n$ stated at the beginning of this section. We prove, using the last two
bounds on $c_n$ and $d_n$, that 
$$
c_n\ll_{\de}(x_c^{-1}+\de)^n,\ n=1,2,\dots,\ \de>0, 
$$
(we use $\ll$ as synonymous to $O(\cdot)$, not to $o(\cdot)$). The simple bounds on $s_n$ and $b_n$ then imply
\begin{eqnarray*}
s_n&=&3b_n\ll_{\de}3\left((x_c^{-1}+\de)^{n+1}+(n+3)(x_c^{-1}+\de)^{n+4}\right)\\
&\ll_{\de}&(x_c^{-1}+\de)^n=(2\cos(\pi/8)+\de)^n,\ \mbox{ for any $\de>0\;$.}
\end{eqnarray*}
To upperbound $c_n$, we fix a $\de\in(0,1)$ and use induction on $n$. We take an $n_0>100$ such that $n>n_0$ 
implies that $(x_c^{-1}+1)^4(8n^2+16n)<(1+\de x_c)^{\sqrt{n/6}}$ and then a constant $\kappa=\kappa(\de)>1$ such 
that $c_n<\kappa$ for $n=1,2,\dots,n_0$. We claim that for every $n=1,2,\dots$,
$$
c_n<\kappa(x_c^{-1}+\de)^n\;.
$$
For $n\le n_0$ this is true. For $n>n_0$ we use the recurrent bound on $c_n$, the upper bound on $d_n$, the 
inductive assumption for $c_{n+2-i}$ and $c_{n+4-i}$ and the property of $n_0$ and get that
\begin{eqnarray*}
c_n&\le&\sum_{i=\sqrt{n/6}}^{n+1} d_ic_{n+2-i}+\sum_{i=\sqrt{n/6}}^{n+3} d_ic_{n+4-i}\\
&<&\kappa\sum_{i=\sqrt{n/6}}^{n+1} 4ix_c^{-i}(x_c^{-1}+\de)^{n+2-i}+
\kappa\sum_{i=\sqrt{n/6}}^{n+3} 4ix_c^{-i}(x_c^{-1}+\de)^{n+4-i}\\
&<&\kappa (x_c^{-1}+1)^4(8n^2+16n)(1+\de x_c)^{-\sqrt{n/6}}(x_c^{-1}+\de)^n\\
&<&\kappa(x_c^{-1}+\de)^n
\end{eqnarray*}
as well. This concludes the proof of the upper bound on $s_n$ and on $\mu$. Theorem~\ref{dcsthm}
is proved.

\section{Proofs of Propositions~\ref{konstrot} and \ref{zjednodus} on winding number}

{\bf Proof of part 1 of Proposition~\ref{konstrot}. }Let $w=(v_1,e_1,\dots,e_n,v_{n+1}=v_1)$ be a closed SAW in the hexagonal graph $H$. 
The \emph{width of $w$} is $s=a_2-a_1$ where $a_2$ (resp. $a_1$) is the maximum (resp. minimum) $x$-coordinate of a vertex in $w$. A \emph{bulge in $w$} is any SAW formed by an interval of edges $e_i,e_{i+1},\dots,e_j$ of $w$ with 
$1\le i<j\le 2n-1, 2\le j-i\le n-1$ (indices $>n$ are understood as reduced modulo $n$) and such that $e_{i+1},e_{i+2},\dots,e_{j-1}$ are all vertical, lying on the line $x=a$, and $e_i,e_j$ are both horizontal, lying to the right of $x=a$. 
Note that $w\cap(x=a_1)$ consists of several (and at least one) disjoint intervals of edges, each of which creates a 
bulge. Hence $w$ contains at least one bulge.

We proceed by induction on $n=|w|\ge6$. If $s=1$ then $w$ is the perimeter of a rectangle $R$ with unit width and 
height $2h$. Traversing the perimeter of $R$ counter-clockwise we see that the vertices on each vertical side 
determine $h+2$ left and $h-1$ right hooks. Thus the weight of all hooks on the perimeter is $2(h+2-(h-1))=6$.  
If the perimeter of $R$ is traversed clockwise, its hooks have weight $-6$. So in this case $r^*(w)=\pm 6$, with 
sign determined by the orientation of $w$.

Let $s>1$. We consider a bulge $B=\{e_i,e_{i+1},\dots,e_j\}$ in $w$ with the maximum value of $a$. None of the horizontal 
edges of $H$ between $e_i$ and $e_j$ is in $w$. Let $I$ be the interval of vertical edges of $H$ on the line $x=a+1$ 
between $e_i$ and $e_j$. The edges in $I$ and not in $w$ form a single nonempty interval $J\subset I$. Indeed, the whole $I$ 
cannot be in $w$ (then we would have $s=1$) and there cannot be two such (neighboring) intervals $J$, because 
their separating interval of edges in $I\cap E(w)$ would create a bulge contradicting the maximality of $a$. We obtain a shorter closed SAW
$w'$ by deleting from $w$ the edges in $B\cup(I\cap E(w))$ and adding the edges in $J$ (oriented compatibly with the 
rest of $w$). Since $w'$ is $2$-regular, connected and self-avoiding, it is indeed a closed SAW. Clearly, $|w'|\le|w|-2$. 
To conclude the argument by induction on length we need to verify that $w$ and $w'$ have the same orientation and 
$r^*(w)=r^*(w')$. 

The exchange of edges does not alter orientation because the rightmost vertical edges of $w$ on $x=a_2$ 
lie to the right of $x=a+1$ and are not affected by the exchange, keep their orientation, and the exterior of $w$ and 
$w'$ lies to the right of them. We calculate the change in $r^*$. Let $w_B$ be the rectangular closed SAW with edges
$B\cup I$ (with orientation inherited from $w$) and $U$ be its vertex set. We have partition $U=U_1\cup\{b,t\}\cup U_2$ 
where $U_1$ are the inner vertices of $J$, $b$ and $t$ are the endvertices of $J$ and $U_2$ are the remaining vertices.
For any $v\in U$ we let $h(v)$ be the hook in $w_B$ with the central vertex $v$. Subscripts $w$ and $w'$, for example 
in $h(b)_w$, define hooks with the given central vertex and contained in the subscripted SAW.
We have
\begin{eqnarray*}
r^*(w)&=&r^*(w')-r(h(b)_{w'})-r(h(t)_{w'})-\sum_{v\in U_1}r(h(v)_{w'})+\\
&&+\;r(h(b)_w)+r(h(t)_w)+\sum_{v\in U_2}r(h(v)_w)\\
&=&r^*(w')-r(h(b)_{w'})-r(h(t)_{w'})+r(h(b)_w)+r(h(t)_w)+\\
&&+\;\sum_{v\in U_1\cup U_2}r(h(v))\\
&=&r^*(w')+\sum_{v\in\{b,t\}}(r(h(v)_w)-r(h(v)_{w'})-r(h(v)))+\sum_{v\in U}r(h(v))\\
&=&r^*(w')+\sum_{v\in\{b,t\}}(r(h(v)_w)-r(h(v)_{w'})-r(h(v)))\pm 6,
\end{eqnarray*}
with the sign $+$ iff $w_B$ is counter-clockwise oriented, because the last but one sum is $r^*(w_B)$ and was 
handled in the case $s=1$. In the last sum we have $r(h(b))=r(h(t))$ and the same for terms with subscripts because 
in the pair of hooks in question one is obtained from the other by mirroring it along the bisector of the segment 
$J$ and then reverting its direction. The definition of $w'$ from $w$ via exchange of edges shows that 
$r(h(v)_{w'})=-r(h(v)_w)$---the edge of $N(v)$ common to both hooks has in them the same orientation. In the hooks 
$h(v)_{w'}$ and $h(v)$ the common edge (which is an extremal edge of $J$) is oriented differently, one hook enters by it 
and the other leaves, and so $r(h(v)_{w'})=r(h(v))$. Finally, the sign of $r(h(t))$, say, equals to that of the last
term since $h(t)$ turns left iff $w_B$ is counter-clockwise oriented. Altogether this shows that the last sum equals 
$\mp 6$ and cancels with the last term. Thus $r^*(w)=r^*(w')$.

{\bf Proof of part 2 of Proposition~\ref{konstrot}. }Let $X\subset V$ be a simple domain, $a,e\in\partial X$ 
be two boundary edges and 
$w_1,w_2\in X(a,e)$ be two SAWs (possibly closed). We distinguish
two cases according to the number of components of $H\backslash X$ the edges $a$ and $e$ lie in, by which we 
mean more precisely their endvertices  in $V\backslash X$. The first case is when $a$ and $e$ lie in the same 
connected component $K$ of $H\backslash X$. Then there is a SAW $w$ going 
from $e$ to $a$ in $K$. The unions $w\cup w_1$ and $w\cup w_2$ are closed SAWs which are moreover oriented in the 
same sense. Suppose for contrary that they have different orientations, one clockwise and the other counter-clockwise.
Then, by Proposition~\ref{jordan}, each vertex in $K$ but not on $w$ is in the interior of $w\cup w_1$ or $w\cup w_2$. 
It follows that for any vertex $v$ in $K$ no matter how we walk from $v$ to infinity the 
walk has to intersect $w\cup w_1\cup w_2$. But this means that $K$ is bounded and hence finite, which was excluded.
Thus, by part $1$, $r^*(w\cup w_1)=r^*(w\cup w_2)=\pm 6$ and from $r^*(w\cup w_1)=r(w)+r(w_1)$ and 
$r^*(w\cup w_2)=r(w)+r(w_2)$ we get that $r(w_1)=r(w_2)$.

The second case is when $a$ is in $K$ and $e$ is in $L$ where $K$ and $L$ are two distinct connected components of 
$H\backslash X$. We will be done as in the previous case if we find a SAW $w$ in $H$ going from $e$ to $a$ such
that it does not intersect, except the overlaps in $e$ and $a$, the SAWs $w_1$ and $w_2$ 
and the resulting closed SAWs $w\cup w_1$ and $w\cup w_2$ have the same orientation. It is not hard to obtain such $w$ using that 
$K$ and $L$ are infinite and thus unbounded. Either there exist vertices in $K$ with $x$-coordinate going to $+\infty$ 
or to $-\infty$ or the same holds with $y$-coordinate, and similarly for $L$. Let us consider the case that both $K$ 
and $L$ contain vertices with $y$-coordinate going to $-\infty$; the other $15$ cases are similar. We fix $c\in\Z$
so small that both $w_1$ and $w_2$ lie above the line $y=c$. Then we walk from $a$ in $K$ down by a SAW $w_3$ until we hit
the line $y=c$ and do the same in $L$ with a SAW $w_4$ starting in $e$. Except the final vertices both $w_3$ 
and $w_4$ lie above $y=c$ and do not intersect as they lie in distinct components of $H\backslash X$. We walk from the 
final vertex of $w_4$ to that of $w_3$ by a SAW $w_5$ that lies on or below $y=c$ and hence besides these final points 
$w_5$ is disjoint to $w_1\cup w_2\cup w_3\cup w_4$. Thus $w=w_3^r\cup w_4\cup w_5$ is the required SAW. If we close 
it by $w_1$ or $w_2$, the orientation is the same because below the lowest horizontal edges of $w_5$ always lies 
the exterior of $w\cup w_1$ and $w\cup w_2$.

{\bf Proof of Proposition~\ref{zjednodus}. }Let $X$ be a simple domain, $a\in\partial X$, $e,f\in A$ be 
two distinct border edges in one part $A$ of a straight partition of $\partial X\backslash\{a\}$ and $U$ be the 
corresponding straight set of vertices of $X$. We may assume that $U$ is a quasisegment joining the vertices of $e$ 
and $f$ in $X$; $e$ and $f$ lie on the same side of $U$. We make $U$ a SAW by orienting it from $e$ to $f$.
Then left and right hooks alternate on $eUf$ except at the ends where the first two and last two hooks have, all three or four, the same orientation. 
It follows from this that in the natural orientation $\ell(\overline{e})=\zeta^3\ell(\overline{f})=-\ell(\overline{f})$. Hence $e$ and $f$ have the same direction when oriented away from $X$. 

We show that $r(e)=r(f)$. Let $w_0$ be the shortest SAW joining $a$ to $U$ in $X$. Thus $w_0$ is disjoint 
to $U$ except that the final vertex $v\in U$. We set $w_1=w_0U'e$ and $w_2=w_0U''f$ where $U'$ is the part 
of $U$ between $v$ and the endvertex of $e$ and similarly for $U''$. So $w_1\in X(a,e)$ and $w_2\in X(a,f)$. Let $h_1$
(resp. $h_2$) be the hook in $w_1$ (resp. $w_2$) with the central vertex $v$. If $v$ is an extremal
vertex of $U$, say an endvertex of $e$, we have $w_1=w_0e$, $w_2=w_0Uf$ and $r(w_2)=r(w_1)-r(h_1)+(-r(h_1))+r(U)+r(h_1)
=r(w_1)+r(U)-r(h_1)=r(w_1)$ because $U$ has an odd number of hooks which alternate on it between left and right and $h_1$
has the same orientation as the two extremal hooks of $U$. If $v$ is an inner vertex of $U$, 
$U=(U')^rU''$ and, by the definition of straight partition, $f(w_0)$ lies on the other side of $U$ than $e$ and $f$. 
Thus, denoting by $h$ the hook in $U$ with the central vertex $v$, we see that again
$r(w_2)=r(w_1)-r(h_1)-r(U'e)+r(h_2)+r(U''f)=r(w_1)+r(e(U')^r)+r(U''f)+r(h)-r(h)-r(h_1)+r(h_2)
=r(w_1)+r(eUf)-r(h)-r(h_1)+r(h_2)=r(w_1)+r(eUf)-3r(h)=r(w_1)$. So $r(w_1)=r(w_2)$ and, by 2 of 
Proposition~\ref{konstrot}, $r(e)=r(f)$.
\eproof

\end{document}